\documentclass[12pt]{amsart} 
\usepackage{latexsym}
\usepackage{amssymb} 
\usepackage{mathrsfs}
\usepackage{amsmath}
\usepackage{latexsym}
\usepackage{delarray} 

\numberwithin{equation}{section}
\theoremstyle{plain}
\newtheorem{thm}{Theorem}[section]
\newtheorem{lem}[thm]{Lemma}
\newtheorem{prop}[thm]{Proposition}
\newtheorem{cor}[thm]{Corollary}
\newtheorem{exam}[thm]{Example}

\newcommand\R{{\mathbb R}}
\newcommand\Rn{{\mathbb R}^n}

\title[A necessary and sufficient condition on scattering]{A necessary 
and sufficient condition on scattering 
for the regularly hyperbolic systems}
\author[Matsuyama and Ruzhansky]
{Tokio Matsuyama${}^\dagger$ and Michael Ruzhansky${}^*$} 
\address{ 
${}^\dagger$Department of Mathematics \endgraf 
Tokai University \endgraf 
Hiratsuka \endgraf 
Kanagawa 259-1292 \endgraf 
Japan \endgraf 
{} \endgraf 
\bigskip
${}^*$Department of Mathematics \endgraf
Imperial College London\endgraf
180 Queen's gate \endgraf 
London SW7 2AZ \endgraf 
United Kingdom}
\thanks{ \endgraf 
The first author was supported by 
Grant-in-Aid for Scientific 
Research (C) (No. 21540198), 
Japan Society for the Promotion of Science. 
The second author was supported in parts by the
EPSRC grant EP/E062873/1 and EPSRC Leadership Fellowship EP/G007233/1.
} 
\email{tokio@keyaki.cc.u-tokai.ac.jp \endgraf 
m.ruzhansky@imperial.ac.uk} 
\begin{document} 

\begin{abstract} 
The present paper is devoted to finding a necessary and sufficient 
condition on the occurence of scattering for the regularly 
hyperbolic systems with time-dependent coefficients whose 
time-derivatives are in $L^1(\mathbb{R})$. More precisely, 
it will be shown that the solutions are asymptotically free if 
the coefficients are stable in the sense of the Riemann integrability 
on $\R$ (R-stability) as $t \to \pm\infty$, while each nontrivial 
solution is never asymptotically free provided that the coefficients 
are not R-stable as $t \to \pm\infty$. As a by-product, the scattering 
operator can be constructed. 
It is expected that the results obtained in 
the present paper would be brought into the study of the asymptotic 
behaviour of Kirchhoff systems. 
\end{abstract}

\maketitle



\section{Introduction.}
\setcounter{equation}{0} 
In this paper we shall give some results on the 
asymptotic behaviour for 
the Cauchy problem of the regularly hyperbolic 
systems with time-dependent coefficients. 
These results would provide a good 
information on the study of the asymptotic 
behaviours for Kirchhoff systems.
In \cite{Matsuyama0} the first author gave the sufficient 
condition on the existence of scattering states for the wave equations, 
and found the special data for the nonexistence 
of scattering states.  
More precisely, there exists a solution $u=u(t,x)$ 
of the Cauchy problem to the strictly hyperbolic equation of second order 
of the form
$$\partial^{2}_t u-c(t)^2 \Delta u=0
$$ 
such that $u$ is never asymptotically free, where we assume that $c(t)$ 
satisfies 
\[
\inf_{t \in \mathbb{R}}c(t)>0, \quad c^\prime (t)\in L^1(\mathbb{R}), 
\quad \lim_{t\to \pm\infty}c(t)=c_{\pm\infty}>0,
\]
\[
\text{$c(t)-c_{\pm\infty}$ is not integrable on $(0,+\infty)$ 
($(-\infty,0)$ resp.).} 
\]
On the contrary, if $c(t)$ is stable, i.e., 
$c(t)-c_{\pm\infty}$ is integrable on $(0,+\infty)$ 
($(-\infty,0)$ resp.), then any solution $u$ is asymptotically free. 
As for the strictly hyperbolic equations of second order 
for ``bounded domains," the similar result was obtained in \cite{Arosio}. 
It should be noted that the results of 
\cite{Matsuyama0} are applied to deduce the nonexistence of scattering 
states for the Kirchhoff equation (see \cite{Matsuyama22}). 
In this sense the behaviour of $c(t)-c_{\pm\infty}$ affects the development 
of the scattering theory for wave equations with time-dependent 
coefficients as well as for the Kirchhoff equation. 

The first order systems often appear in the analysis of 
equations of orders larger than two and in the analysis 
of coupled equations of second order (see examples 
\ref{exam:Example 1}--\ref{exam:Example 2} below). In the present paper 
we will find the necessary 
and sufficient condition on the occurrence 
of scattering for the regularly 
hyperbolic system with time-dependent coefficients, 
which extends the results of \cite{Matsuyama0}. 
We will also construct the wave operators and scattering operators 
by using the asymptotic integrations method, which were developed 
in \cite{MR2}. In hyperbolic systems we will impose 
a stability condition on the characteristic roots of the symbol 
of the differential operator. 

Apart from the scattering problem, 
the dispersion for hyperbolic systems is also of great interest. 
Large hyperbolic systems appear in many applications, for example
the Grad systems of gas dynamics, hyperbolic systems
in the Hermite-Grad decomposition
of the Fokker-Planck equation, etc. Thus,
for general hyperbolic equations and systems with constant
coefficients a comprehensive analysis of dispersive and Strichartz
estimates was carried out in \cite{RS}. The dispersion for
scalar equations based on the asymptotic integration method
was analysed by the authors in \cite{MR2}, motivated by the
higher order Kirchhoff equations.
The dispersion for hyperbolic systems with time-dependent
coefficients will be discussed in \cite{MR3} and 
will appear elsewhere, as well as the applications of the
obtained results to the Kirchhoff systems.

To be more precise, 
let us consider the Cauchy problem 
\begin{equation}
D_t U=A(t,D_x)U \qquad \text{with $D_t=-i\partial_t$ and 
$D_{x_j}=-i\partial_{x_j}$ $(j=1,\ldots,n)$,}
\label{Eq}
\end{equation}
$i=\sqrt{-1}$, for $t\ne0$, with Cauchy data 
\begin{equation}
U(0,x)={}^T (f_0(x),\ldots,f_{m-1}(x))\in (L^2(\mathbb{R}^n))^m. 
\label{Cauchy data}
\end{equation} 
The operator $A(t,D_x)$ is the first order $m\times m$ 
pseudo-differential\footnote{We note that it is 
important to allow $a_{ij}$ to be 
pseudo-differential here since we want the results to hold for 
scalar higher order equations as well, e.g. see
Example \ref{exam:Example 2}.} system,
with symbol $A(t,\xi)$ of the form
$A(t,\xi)=\{a_{ij}(t,\xi)\}_{i,j=1}^n$, where we assume that 
$a_{ij}(t,\xi)$ are positively homogeneous of order one in $\xi$, 
$a_{ij}(t,\lambda \xi)=\lambda a_{ij}(t,\xi)$ for 
$\lambda>0$, $\xi\in\Rn\backslash 0$,
and satisfy 
\begin{equation}\label{hyp1}
\text{$a_{ij}(\cdot,\xi) \in \mathrm{Lip}_{\mathrm{loc}}(\R)$,
and $\partial_t a_{ij}(\cdot,\xi) \in L^1(\R)$ for $i,j=1,\ldots,m$.}
\end{equation}
We will also assume that the symbol $A(t,\xi)$ satisfies the regularly 
hyperbolic condition in the sense of Mizohata (\cite{Mizohata}): 
\begin{equation}
\text{$\mathrm{det} (\tau I-A(t,\xi))=0$ has (in $\tau$) 
real and distinct roots 
$\varphi_1(t,\xi),\ldots,\varphi_m(t,\xi)$,}
\label{hyp2}
\end{equation}
i.e., 
\begin{equation}\label{hyp3}
\underset{t\in \mathbb{R},|\xi|=1,j\ne k}{\inf}
|\varphi_j(t,\xi)-\varphi_k(t,\xi)|=d>0. 
\end{equation}
Notice that each characteristic root $\varphi_j(t,\xi)$ is 
positively homogeneous of order one in $\xi$. 
The assumption \eqref{hyp1} assures the existence of the 
limiting functions 
$a_{ij}^\pm(\xi)$, $i,j=1,\ldots,m$, such that 
\begin{equation}\label{EQ:asymptotic}
a_{ij}(t,\xi)\to a_{ij}^\pm(\xi) \quad (t\to \pm\infty),
\end{equation}
and we can expect that the solution $U(t,x)$ of 
\eqref{Eq}--\eqref{Cauchy data} is 
asymptotic to some solution 
of the following hyperbolic system with constant coefficients 
as $t \to \pm\infty$: 
\begin{equation}
D_t V=A_\pm(D_x)V,  
\label{Limiting Eq}
\end{equation}
where $A_\pm(D_x)$ is the $m\times m$ first order 
pseudo-differential system, with symbol 
$A_\pm(\xi)=\{a_{ij}^\pm(\xi)\}_{i,j=1}^m$. 
Since the characteristic roots depend continuously 
on the coefficients, it 
follows from \eqref{hyp1}--\eqref{EQ:asymptotic} 
that the operator $D_t-A_\pm(D_x)$ 
is regularly hyperbolic. This is because it will be shown in 
Proposition \ref{prop:root} 
below that there exists the limiting phases $\varphi_j^\pm(\xi)$ 
of $\varphi_j(t,\xi)$ for $j=1,\ldots,m$:
\[
\lim_{t\to \pm\infty}\varphi_j(t,\xi)=\varphi_j^\pm(\xi), 
\quad j=1,\ldots,m.
\] 
Hence by using \eqref{hyp3}, we have also 
\begin{equation}\label{Limiting phase}
\underset{|\xi|=1,j\ne k}{\inf}
|\varphi_j^\pm(\xi)-\varphi_k^\pm(\xi)|=d>0.
\end{equation}

We are now in a position to state our results. For this purpose, we recall 
the notion of scattering states. 
We say that the solution $U(t,x)$ of 
$D_tU=A(t,D_x)U$ is asymptotically free in 
$(L^2(\mathbb{R}^n))^m$, 
if it is asymptotic to some solution $V_\pm(t,x)$ of 
$D_t V=A_\pm(D_x)V$, i.e., 
\[
\| U(t,\cdot)-V_\pm(t,\cdot) \|_{(L^2(\mathbb{R}^n))^m} \to 0 
\quad (t \to \pm\infty). 
\]

We shall prove here the following theorem which gives a necessary 
and sufficient condition on the existence of scattering states. 
  
\begin{thm} \label{thm:Nonfree}
Assume \eqref{hyp1}--\eqref{hyp3}. 
Then the integrabililty condition 
on $\varphi_j(t,\xi)-\varphi^{\pm}_j(\xi)$ is necessary and sufficient for the 
asymptotically free property for 
\eqref{Eq}--\eqref{Cauchy data}. More precisely, the following 
assertions hold{\rm :} 

{\rm i)} If $\varphi_j(t,\xi)-\varphi^{\pm}_j(\xi)$ $(j=1,\ldots,m)$ 
are integrable on 
$(0,+\infty)$ $\textrm{$($resp. $(-\infty,0))$,}$ i.e., the functions 
\begin{equation}\label{EQ:cond-ns}
\psi_{j,\pm}(t,\xi)\equiv 
\int^{t}_0 \left( \varphi_j(s,\xi)-\varphi^{\pm}_j(\xi) \right) 
\, ds, \quad j=1,\ldots,m, 
\end{equation} 
have the finite limits for each $\xi \ne 0$ as $t\to\pm\infty$, 
then each solution 
$U(t,x)\in C(\R;(L^2(\mathbb{R}^n))^m)$ 
of \eqref{Eq}--\eqref{Cauchy data} is 
asymptotically free in $(L^2(\mathbb{R}^n))^m$. 

Moreover, the mapping $($the inverse of the wave operators 
$\mathscr{W}_\pm)$ 
\[
\mathscr{W}_\pm^{-1}:U(0) \mapsto V_\pm(0)
\]
is well-defined and bounded on $(L^2(\Rn))^m$.

{\rm ii)} If $\psi_{j,\pm}(t,\xi)$ satisfy 
\begin{equation}
\lim_{t\to \pm\infty} \vert \psi_{j,\pm}(t,\xi)
\vert=+\infty, \quad \xi \ne 0, 
\quad \quad j=1,\ldots,m,
\label{nonfree}
\end{equation}
then each non-trivial solution 
$U(t,x)\in C(\R;(L^2(\mathbb{R}^n))^m)$ of 
\eqref{Eq}--\eqref{Cauchy data} is never asymptotically free in 
$(L^2(\mathbb{R}^n))^m$. 
\end{thm} 
We note that condition \eqref{EQ:cond-ns} is stronger than
\eqref{EQ:asymptotic}, thus indeed specifying
a subset of systems. There are different 
sufficient criteria for
\eqref{EQ:cond-ns} to hold. For example, for
$A(\cdot,\xi)\in C^1(\R)$, if for all $|\xi|=1$ we have 
$$t(\psi_j(t,\xi)-\psi_j^{\pm}(\xi))=o(1) \textrm{ as } t\to\pm\infty,
\textrm{ and }
t\psi_j'(t,\xi)\in L^1(\R),\quad j=1,\ldots,m,$$
then \eqref{EQ:cond-ns} follows. Indeed, in this case
we have $\psi_j(\cdot,\xi)\in C^1(\R)$ for all $j$,
and the statement follows from the trivial identity
$$
\int_0^t s \psi_j'(s,\xi)\ ds=
\int_0^t (\psi_j^\pm(\xi)-\psi_j(s,\xi))\ ds+
t(\psi_j(t,\xi)-\psi_j^\pm(\xi)).
$$
Thus, under the assumption $t(\psi_j(t,\xi)-\psi_j^{\pm}(\xi))=o(1)$
as $t\to\pm\infty$,
condition \eqref{EQ:cond-ns} is equivalent to 
$t\psi_j'(t,\xi)\in L^1(\R), j=1,\ldots,m$.

Next we state the results on the existence of wave operators. 
Let us consider 
the Cauchy problem for the regularly hyperbolic system with 
constant coefficients 
\begin{equation}\label{inverse Equation}
D_t V_\pm=A_\pm(D_x)V_\pm, \quad x\in \Rn, \quad \pm t> 0,
\end{equation}
with Cauchy data 
\begin{equation}\label{EQ:pm-Cauchy}
V_\pm(0,x)={}^T (f_0^\pm(x),\ldots,f_{m-1}^\pm(x)),
\end{equation}
where $A_\pm(D_x)$ is the pseudo-differential operator with 
symbol $\{a_{ij}^\pm(\xi)\}_{i,j=1}^m$. We will assume that 
the characteristic roots $\varphi_1^\pm(\xi),\ldots,\varphi_m^\pm(\xi)$ 
of the operator $D_t-A_\pm(D_x)$ are real and distinct, 
i.e., 
\begin{equation}
\mathrm{det}(\tau I-A_\pm(\xi))
=(\tau-\varphi_1^\pm(\xi)) \cdots (\tau-\varphi_m^\pm(\xi)), 
\label{strict hyperbolicity5}
\end{equation} 
\begin{equation}
\inf_{|\xi|=1,j \ne k} |\varphi_j^\pm(\xi)-\varphi_k^\pm(\xi)|=d>0.
\label{strict hyperbolicity6}
\end{equation}
Then the following theorem assures the existence 
of wave operators. 
\begin{thm}\label{thm:Scattering}
Assume \eqref{strict hyperbolicity5}--\eqref{strict hyperbolicity6}. 
Suppose that $a_{ij}(t,\xi)$ are positively
homogeneous of order one in $\xi$, and 
satisfy \eqref{hyp1} 
in such a way that 
\[
a_{ij}(t,\xi) \to a_{ij}^\pm(\xi) \quad \text{for all $i,j$ 
$(t \to \pm\infty)$.}
\]
Let $\varphi_1(t,\xi),\ldots,\varphi_m(t,\xi)$ be the characteristic 
roots of the regularly hyperbolic 
operator $D_t-A(t,D_x)$, with symbol $A(t,\xi)=\{a_{ij}(t,\xi)\}_{i,j=1}^m$. 
Assume that 
\[
\text{$\varphi_j(t,\xi)-\varphi_j^\pm(\xi)$ is integrable on 
$(0,+\infty)$ $($resp. $(-\infty,0))$ for each $\xi\ne0$.}
\]
Then for any solution $V_\pm(t,x)\in C(\R;(L^2(\mathbb{R}^n))^m)$ of 
\eqref{inverse Equation}--\eqref{EQ:pm-Cauchy}, 
there exists a unique solution 
$U(t,x)\in C(\R;(L^2(\mathbb{R}^n))^m)$ of 
$D_t U=A(t,D_x)U$ such that 
\[
\| V_\pm(t)-U(t)\|_{(L^2(\mathbb{R}^n))^m} \to 0 
\quad (t \to \pm\infty).
\] 

Moreover, the mapping $($wave operator$)$ 
\[
\mathscr{W}_\pm: V_\pm(0) \mapsto U(0)
\]
is well-defined and bounded on $(L^2(\Rn))^m$.
\end{thm}

As a consequence of Theorems \ref{thm:Nonfree}--\ref{thm:Scattering}, 
we can construct the scattering operators. More precisely, we have: 
\begin{cor}
Assume that $a_{ij}^+(\xi)=a_{ij}^-(\xi)$ hold for 
$i,j=1,\ldots,m$ in Theorems 
\ref{thm:Nonfree}--\ref{thm:Scattering}. 
Then the mapping 
\[
S=\mathscr{W}_+^{-1}\mathscr{W}_- : V_-(0)\mapsto V_+(0)
\] 
defines the scattering operator, and it 
is bijective and bounded on $(L^2(\mathbb{R}^n))^m$. 
\end{cor}

Finally, let us look at some examples to which our theorems can be applied.
We note that although the equations may be of high order, it is important
that we impose conditions only on one time-derivative of the coefficients.
This is of crucial importance to being able to apply the obtained
results to the Kirchhoff equations.

Our first example deals with higher order scalar equations.
\begin{exam}\label{exam:Example 1}
Consider the Cauchy problem to the $m^{\mathrm{th}}$ order 
strictly hyperbolic 
equation 
\[
L(t,D_t,D_x)u 
\equiv D^{m}_t u+\sum_{\underset{j \le m-1}{\vert \nu \vert+j=m}} 
a_{\nu,j}(t) D^{\nu}_x D^{j}_t u=0, \quad t \ne 0,
\]
with Cauchy data 
\[
D^{k}_t u(x,0)=f_k(x) \in H^{m-1-k}(\Rn), \quad 
k=0,1,\cdots,m-1, \quad x \in \mathbb{R}^n, 
\]
where $D_t=-i\partial_t$ and $D^{\nu}_x=
\left(-i\partial_{x_1} \right)^{\nu_1}
\cdots \left(-i\partial_{x_n} \right)^{\nu_n}$, 
$i=\sqrt{-1}$, for $\nu=(\nu_1,\ldots,\nu_n)$. We assume that 
$a_{\nu,j}(t)$ belong to $\mathrm{Lip}_{\mathrm{loc}}(\mathbb{R})$ 
and satisfy 
\[
a_{\nu,j}^\prime(t) \in L^1(\mathbb{R}) \quad \text{for all $\nu,j$, 
and $k=1,\ldots,m-1$},
\]
and the symbol $L(t,\tau,\xi)$ of the operator 
$L(t,D_t,D_x)$ has real roots 
$\varphi_1(t,\xi),\ldots,\varphi_m(t,\xi)$ which are
uniformly distinct for $\xi \ne0$, i.e., 
\[
L(\tau,\xi)=(\tau-\varphi_1(t,\xi)) \cdots (\tau-\varphi_m(t,\xi)), 
\]
\[
\inf_{\underset{j \ne k}{\vert \xi \vert=1,t \in \mathbb{R}}} 
\vert \varphi_j(t,\xi)-\varphi_k(t,\xi) \vert=d>0.
\]
The reference equation is 
\[
D^{m}_t v_\pm+\sum_{\underset{j \le m-1}{\vert \nu \vert+j=m}} 
a_{\nu,j}^\pm D^{\nu}_x D^{j}_t v_\pm=0,
\]
with $a_{\nu,j}^\pm=\lim_{t\to\pm\infty}a_{\nu,j}(t)$, and the energy space 
is $\dot{H}^{m-1}(\Rn)\times \cdots \times L^2(\Rn)$. 
\end{exam}

The following example deals with coupled second order equations.
\begin{exam}\label{exam:Example 2}
Let us consider the Cauchy problem 
\[
\partial_t^2u-c_1(t)^2\Delta u+P_1(t,D_x)v=0,
\]
\[
\partial_t^2v-c_2(t)^2\Delta v+P_2(t,D_x)u=0,
\]
for some second order 
homogeneous polynomials $P_1(t,D_x),P_2(t,D_x)$
which may depend on time, 
where we assume that 
\[
c_k(t), P_k(t,\xi)\in \mathrm{Lip}_{\mathrm{loc}}(\R), \quad 
c_k^\prime(t), P_k^\prime(t,\xi)\in L^1(\R), \qquad (k=1,2),\quad \xi\in\Rn,
\]
\[
\inf_{t\in \R,|\xi|=1} 
\left((c_1(t)^2-c_2(t)^2)^2+4P_1(t,\xi)P_2(t,\xi)\right)>0, 
\]
\[
\inf_{t\in \R,|\xi|=1} \left(c_1(t)^2c_2(t)^2-P_1(t,\xi)P_2(t,\xi)
\right)>0.
\]
By taking the Fourier transform in the space variables and introducing the vector 
\[
V(t,\xi)={}^T(v_1(t,\xi),v_2(t,\xi),v_3(t,\xi),v_4(t,\xi))
={}^T(|\xi|\widehat{u}(t,\xi),\widehat{u}^\prime(t,\xi),
|\xi|\widehat{v}(t,\xi),\widehat{v}^\prime(t,\xi))
\]
we obtain the system 
\begin{align*}
\frac{\partial V}{\partial t}=& i
\begin{pmatrix}
0&-i|\xi|&0&0\\
ic_1(t)^2|\xi|&0& iP_1(t,\xi)|\xi|^{-1}&0\\
0&0&0& -i|\xi|\\
iP_2(t,\xi)|\xi|^{-1} &0& ic_2(t)^2|\xi|&0
\end{pmatrix}V\\
=& iA(t,\xi)V. 
\end{align*}
The four characteristic roots of $\mathrm{det}(\tau I-A(t,\xi))=0$ in $\tau$
are given by
$$
\varphi_{1,2,3,4}(t,\xi)=
\pm \frac{|\xi|}{\sqrt{2}}
\sqrt{c_1(t)^2+c_2(t)^2\pm
\sqrt{(c_1(t)^2-c_2(t)^2)^2+4P_1(t,\xi) P_2(t,\xi)|\xi|^{-4}}}.
$$ 
The reference system is 
\[
\partial_t^2u_\pm-c_{1,\pm}^2\Delta u_\pm+P_{1,\pm}(D_x)v_\pm=0,
\]
\[
\partial_t^2v_\pm-c_{2,\pm}^2\Delta v_\pm+P_{2,\pm}(D_x)u_\pm=0,
\]
with the limits $c_{k,\pm}=\lim_{t\to\pm\infty}c_k(t)$,
$P_{k,\pm}(\xi)=\lim_{t\to\pm\infty}P_k(t,\xi)$ for $k=1,2$, 
and the energy space is $(\dot{H}^1(\Rn)\times L^2(\Rn))^2$. 
\end{exam}

We conclude this section by stating our plan. In \S 2 we will find the 
representation formulae for \eqref{Eq}--\eqref{Cauchy data}. The 
proof of Theorem \ref{thm:Nonfree} 
will be given in \S 3 and \S 4. 
In the last section we will prove Theorem \ref{thm:Scattering}.

\section{Representation formulae via asymptotic integrations} 
In this section we will derive the 
representation formulae for \eqref{Eq} 
along the method of \cite{MR2}. 
Let us first analyse certain basic properties of characteristic roots 
$\varphi_k(t,\xi)$ of \eqref{hyp2}. The first part of
the following statement was established 
in \cite{MR3}. For the completeness, we will give the proof. 

\begin{prop} \label{prop:root}
Let the operator $D_t-A(t,D_x)$ satisfy the properties
\eqref{hyp2}--\eqref{hyp3}. Then each 
$\partial_t \varphi_k(t,\xi)$, $k=1,\ldots,m$, is 
homogeneous of order one in $\xi$, and there exist a constant $C>0$ such that 
\begin{equation}\label{EQ:phi-bounded0}
| \partial_t\varphi_k(t,\xi)|\leq C|\xi| 
\quad \textrm{for all}\;\; \xi\in\Rn, 
\;\; t \in \mathbb{R}, \;\; k=1,\ldots,m.
\end{equation}
Moreover, if $\partial_t a_{ij}(\cdot,\xi)\in L^1(\R)$ for all 
$\xi\in \Rn$ and $i,j=1,\ldots,m$, 
then we have also 
$\partial_t\varphi_k(\cdot,\xi)\in L^1(\R)$ for all 
$\xi\in\mathbb{R}^n$. Furthermore, there exist functions 
$\varphi_k^\pm\in C^\infty(\Rn\backslash 0)$, homogeneous of order one, 
such that 
\begin{equation}\label{EQ:convphi}
\varphi_k(t,\xi)\to 
\varphi_k^\pm(\xi) \qquad \textrm{ as } 
t\to\pm\infty,
\end{equation} 
for all $\xi\in\Rn$, and $k=1,\ldots,m$. 
Finally, we have the following 
formula for the derivatives of characteristic roots{\rm :}
\begin{equation}
\partial_t \varphi_k(t,\xi)= 
-\sum_{j=0}^m \partial_t \alpha_{m-j}(t,\xi) \varphi_k(t,\xi)^j
\prod_{r\not=k} 
\left( \varphi_k(t,\xi)-\varphi_r(t,\xi)\right)^{-1},
\label{1-product0}
\end{equation}
where 
\begin{eqnarray*}
\alpha_k(t,\xi)=(-1)^k \sum_{i_1<i_2<\cdots<i_k}
\mathrm{det}\left(
\begin{array}{ccc}
\displaystyle{\sum_{r=1}^n} a_{i_1 i_1}(t,\xi) 
& \cdots & \displaystyle{\sum_{r=1}^n}
a_{i_1 i_k}(t,\xi) \\
\vdots & \ddots & \vdots \\
\displaystyle{\sum_{r=1}^n} a_{i_k i_1}(t,\xi)
& \cdots & \displaystyle{\sum_{r=1}^n}
 a_{i_k i_k}(t,\xi) \\
\end{array}
\right).
\end{eqnarray*}
\end{prop}
\begin{proof}
Let us show first that $\varphi_k(t,\xi)$ is bounded with respect to 
$t\in \mathbb{R}$, i.e., 
\begin{equation}
\vert \varphi_k(t;\xi) \vert \le C \vert \xi \vert, \quad 
\text{for all $\xi \in \mathbb{R}^n$, $t \in \mathbb{R}$, $k=1,\ldots,m$.}
\label{phibdd}
\end{equation}
We will use the fact that $\varphi_k(t,\xi)$ are roots of the polynomial 
$L(t,\tau,\xi)=\mathrm{det}(\tau I-A(t,\xi))$ of the form 
$L(t,\tau,\xi)=\tau^m+\alpha_1(t,\xi)\tau^{m-1}
+\cdots+\alpha_m(t,\xi)$ with $|\alpha_j(t,\xi)|\leq M|\xi|^j$, 
for some $M\geq 1$. 
Suppose that one of its roots $\tau$ satisfies $|\tau(t,\xi)|>2M|\xi|$. 
Then
\begin{align*}|L(t,\tau,\xi)| & 
 \geq |\tau|^{m} \left(1-\frac{|\alpha_1(t,\xi)|}{|\tau|}-
\cdots-\frac{|\alpha_m(t,\xi)|}{|\tau|^m} \right) \\
& \geq 
2M|\xi|^m \left(1-\frac{1}{2}-\frac{1}{4M}-\cdots-\frac{1}{2^m M^{m-1}}
\right)>0,
\end{align*}
hence $|\tau(t,\xi)|\leq 2M|\xi|$ for all $\xi\in\Rn$. 
Thus we establish \eqref{phibdd}. 

Differentiating \eqref{hyp2} with respect to $t$, we get
$$ 
\frac{\partial L(t,\tau,\xi)}{\partial t}=
\sum_{j=0}^m \partial_t \alpha_{m-j}(t,\xi) \tau^j=
-\sum_{k=1}^m \partial_t\varphi_k(t,\xi)\prod_{r\not=k} 
\left(\tau-\varphi_r(t,\xi)\right).
$$
Setting $\tau=\varphi_k(t,\xi)$, we obtain
\begin{equation}
\partial_t\varphi_k(t;\xi)\prod_{r\not=k} 
\left(\varphi_k(t,\xi)-\varphi_r(t,\xi)\right)= -
\sum_{|\nu|+j=m} \partial_t \alpha_{m-j}(t,\xi) \varphi_k(t,\xi)^j,
\label{2-product}
\end{equation}
implying \eqref{1-product0}.
Now, using \eqref{hyp3}, \eqref{phibdd}, 
and the assumption 
that $\partial_t \alpha_j(\cdot,\xi)\in L^1(\R_t)$ for all $j$, 
we conclude that \eqref{EQ:phi-bounded0} holds and 
$\partial_t\varphi_k(\cdot,\xi)\in L^1(\R)$ for all $\xi\in\Rn$ and 
$k=1,\ldots,m$. 
The homogeneity of order one of $\partial_t \varphi_\ell(t,\xi)$ 
is an immediate consequence of \eqref{2-product} and its derivatives. 

Finally, setting $\varphi_k^\pm(\xi)=\varphi_k(0,\xi)
+\int_0^{\pm\infty}\partial_t\varphi_k(t,\xi) \, dt$, 
we get \eqref{EQ:convphi}. 
The proof is complete.
\end{proof}

We prepare the next lemma. 
\begin{lem}[\cite{Mizohata}~Proposition~6.4] \label{Diagonalisation}
Assume \eqref{hyp1}--\eqref{hyp3}. 
Then there exists a matrix $\mathscr{N}=\mathscr{N}(t,\xi)$ of homogeneous 
degree $0$ in $\xi$ satisfying the following properties{\rm :} 

\noindent 
{\rm (i)} $\mathscr{N}(t,\xi) A(t,\xi/|\xi|)
=\mathscr{D}(t,\xi)\mathscr{N}(t,\xi)$, 
where 
\[
\mathscr{D}(t,\xi)
=\mathrm{diag} \left( \varphi_1(t,\xi/|\xi|),\ldots,
\varphi_m(t,\xi/|\xi|)\right);
\]

\noindent 
{\rm (ii)} $\displaystyle
{\inf_{\xi \in \mathbb{R}^n\backslash 0, t \in \mathbb{R}}}
\vert {\rm det} \, \mathscr{N}(t,\xi)) \vert>0;$

\noindent 
{\rm (iii)} $\mathscr{N}(t,\xi)$ belongs to 
$\mathrm{Lip}_{\mathrm{loc}}
\left(\R_t;(C^\infty(\Rn_\xi\setminus0))^{m^2}\right)$ and 
\[\partial_t \mathscr{N}(t,\xi)\in (L^1(\mathbb{R}))^{m^2} 
\quad \text{
for each $\xi\not=0$.}
\]
\end{lem}

Applying the Fourier transform on $\mathbb{R}_x^n$, we get the 
following ordinary differential system from \eqref{Eq}: 
\begin{equation} \label{EQ}
D_t V=A(t,\xi/|\xi|)|\xi|V. 
\end{equation}
We find the asymptotic integration of \eqref{EQ} following 
Ascoli \cite{Ascoli} and Wintner \cite{Wintner}, cf. Hartman \cite{Hartman}. 
Multiplying \eqref{EQ} by $\mathscr{N}=\mathscr{N}(t,\xi)$ 
equations from Lemma \ref{Diagonalisation} and putting $\mathscr{N}V=W$, we get 
\begin{equation}
D_t W=\mathscr{D}|\xi|W+(D_t \mathscr{N})V
=\left(\mathscr{D}|\xi|+(D_t \mathscr{N})\mathscr{N}^{-1}\right)W,
\label{Ref EQ}
\end{equation}
since $\mathscr{N}A(t,\xi/|\xi|)=\mathscr{D}\mathscr{N}$ by 
Lemma \ref{Diagonalisation}. We can expect that the solutions of 
\eqref{Ref EQ} are asymptotic to some solution of 
\begin{equation}
D_t \pmb{y}=\mathscr{D}|\xi|\pmb{y}.
\label{Ref EQ1}
\end{equation}
Let $\Phi(t,\xi)$ be the fundamental matrix of \eqref{Ref EQ1}, i.e., 
\[
\Phi(t,\xi)=\mathrm{diag}\left( 
e^{i\int_0^t \varphi_1(s,\xi) \, ds}, 
\cdots, e^{i\int_0^t \varphi_m(s,\xi) \, ds}\right). 
\]
If we perform the Wronskian transform 
$\pmb{a}(t,\xi)=\Phi(t,\xi)^{-1}W(t,\xi)$, then the 
system \eqref{Ref EQ} reduces to the system 
$D_t\pmb{a}=C(t,\xi)\pmb{a}$, where $C(t,\xi)$ is given by 
\[
C(t,\xi)=\Phi(t,\xi)^{-1} (D_t \mathscr{N}(t,\xi))
\mathscr{N}(t,\xi)^{-1}\Phi(t,\xi). 
\]
We note that $C(\cdot,\xi)\in (L^1(\R))^{m^2}$, since 
$D_t \mathscr{N}(\cdot,\xi)\in (L^1(\R))^{m^2}$ by Lemma \ref{Diagonalisation}.
Hence $D_t\pmb{a}(\cdot,\xi)\in (L^1(\R))^m$; thus there exist the limits 
\[
\lim_{t\to\pm\infty}\pmb{a}(t,\xi)=\pmb{\alpha}_{\pm}(\xi).
\]
Since $W(t,\xi)=\Phi(t,\xi)\pmb{a}(t,\xi)$ and 
$\mathscr{N}(t,\xi)V(t,\xi)=W(t,\xi)$, we get 
\[
V(t,\xi)=\mathscr{N}(t,\xi)^{-1}\Phi(t,\xi)
\pmb{a}(t,\xi). 
\]

Now let $(V_0(t,\xi),\ldots,V_{m-1}(t,\xi))$ be the fundamental 
matrix of \eqref{EQ}. This means, in particular, that 
$(V_0(0,\xi),\ldots,V_{m-1}(0,\xi))=I$. Then each $V_j(t,\xi)$ can be 
represented by 
\[
V_j(t,\xi)=\mathscr{N}(t,\xi)^{-1}\Phi(t,\xi)\pmb{a}^j(t,\xi), 
\]
where $\pmb{a}^j(t,\xi)$ are the corresponding amplitude functions to 
$V_j(t,\xi)$. 
Since $\widehat{U}(t,\xi)=\sum_{j=0}^{m-1} 
V_j(t,\xi)\widehat{f}_j(\xi)$, 
we arrive at 
\[
\widehat{U}(t,\xi)=\sum_{j=0}^{m-1} \mathscr{N}(t,\xi)^{-1}\Phi(t,\xi)
\pmb{a}^j(t,\xi)\widehat{f}_j(\xi).
\]
Finally, let us find the estimates of the amplitude functions 
$\pmb{a}^j(t,\xi)$. Recalling that $\pmb{a}^j(t,\xi)$ satisfy the 
problem 
\[
D_t\pmb{a}^j=C(t,\xi)\pmb{a}^j \quad \text{with 
$(\pmb{a}^0(0,\xi),\cdots,\pmb{a}^{m-1}(0,\xi))=\mathscr{N}(0,\xi)$,}
\] 
we can write $\pmb{a}^j(t,\xi)$ by the Picard series: 
\[
\pmb{a}^j(t,\xi)=\left(I+i\int^t_0 C(\tau_1,\xi) \, d \tau_1
+i^2\int^{t}_0 C(\tau_1,\xi) \, 
d \tau_1\int^{\tau_1}_0 
C(\tau_2,\xi)\, d \tau_2+\cdots \right) \pmb{a}^j(0,\xi).
\]
This implies that 
\[
\left|\pmb{a}^j(t,\xi)\right| \le 
e^{c\int_\R 
\|\partial_t \mathscr{N}(s,\xi)\|_{L^\infty(\Rn)}\, ds} 
|\pmb{a}^j(0,\xi)|,
\]
where we have used the following: 

\smallskip

{\bf Fact.} {\em 
Let $f(t)$ be a continuous function on $\mathbb{R}$. Then 
\[
e^{\int^{t}_{s} f(\tau) \, d \tau}
=1+\int^{t}_{s} f(\tau_1) \, d \tau_1
+\int^{t}_{s} f(\tau_1) \, 
d \tau_1 \int^{\tau_1}_{s} f(\tau_2) 
\, d \tau_2+\cdots.
\]
}

\smallskip

Summarising the above argument, we obtain 
\begin{prop} \label{prop:Rep}
Assume \eqref{hyp1}--\eqref{hyp3}. 
Let 
$\mathscr{N}(t,\xi)$ be the diagonaliser of 
$A(t,\xi/|\xi|)$ constructed in Lemma \ref{Diagonalisation}. 
Then there exist vector-valued functions 
$\pmb{a}^j(t,\xi)$, $j=0,1,\ldots,m-1$, 
determined by the initial value problem 
\[
D_t\pmb{a}^j(t,\xi)=C(t,\xi)\pmb{a}^j(t,\xi), \qquad 
\left(\pmb{a}^1(0,\xi),\cdots,\pmb{a}^m(0,\xi)\right)=\mathscr{N}(0,\xi),
\]
\[
\text{with} \quad C(t,\xi)=\Phi(t,\xi)^{-1} (D_t \mathscr{N}(t,\xi))
\mathscr{N}(t,\xi)^{-1}\Phi(t,\xi) \in (L^1(\mathbb{R}_t))^{m^2},
\]
such that 
the solution $U(t,x)$ of \eqref{Eq} is represented by 
\begin{equation}\label{EQ:Representation of U}
U(t,x)=\sum_{j=0}^{m-1} 
\mathscr{F}^{-1} \left[\mathscr{N}(t,\xi)^{-1}\Phi(t,\xi)
\pmb{a}^j(t,\xi)\widehat{f}_j(\xi) \right](x). 
\end{equation}
Moreover, the limits 
\[
\lim_{t\to \pm\infty}\pmb{a}^j(t,\xi)=\pmb{\alpha}^j_\pm(\xi), 
\quad j=0,1,\ldots,m-1, 
\]
exist, and there exists a constant $c>0$ such that 
\[
\left|\pmb{a}^j(t,\xi)\right| \le c, \quad j=0,1,\ldots,m-1, 
\]
for all $t\in \R$ and $\xi\in \Rn$. 
\end{prop}
Proposition \ref{prop:Rep} is known as Levinson's lemma 
(see Coddington and Levinson \cite{Coddington}) in the 
theory of ordinary differential equations; the new feature
here is the additional dependence on $\xi$, which is crucial
for our analysis.

\section{Proof of Theorem \ref{thm:Nonfree} (i).} 
By our assumption that $\varphi_j(t;\xi)-\varphi^{\pm}_j(\xi)$, 
$j=1,\ldots,m$, are 
integrable on $(0,+\infty)$ ($(-\infty,0)$ resp.), we can define 
functions $\Theta^{\pm}_j(\xi)$ to be 
\[
\Theta^{\pm}_j(\xi)=\int^{\pm\infty}_{0} \left(\varphi_j(s,\xi)
-\varphi^{\pm}_j(\xi) \right) \, ds.
\]
Put 
\[
\Phi_j(s,\xi)=\varphi_j(s,\xi)-\varphi^{\pm}_j(\xi). 
\]
Then we can write 
\[
\int^{t}_{0}\varphi_j(s,\xi) \, ds
=\varphi^{\pm}_j(\xi)t+\Theta^{\pm}_j(\xi)
-\int^{\pm\infty}_{t} \Phi_j(s,\xi)\, ds. 
\]
Therefore, we have 
\begin{align}
e^{i\int^{t}_{0}\varphi_j(s,\xi) \, ds}
&=e^{i\left(\varphi^{\pm}_j(\xi)t+\Theta^{\pm}_j(\xi) \right)}
+e^{i\left(\varphi^{\pm}_j(\xi)t+\Theta^{\pm}_j(\xi)\right)}
\left( \exp \left(-i\int^{\pm\infty}_{t} \Phi_j(s,\xi)\, ds 
\right)-1\right) \nonumber \\ 
&\equiv e^{i\left(\varphi^{\pm}_j(\xi)t+\Theta^{\pm}_j(\xi) 
\right)}+\Psi_j(t,\xi), \label{phase}
\end{align}
with 
\[
\Psi_j(t,\xi)=O\left(\int^{+\infty}_{\vert t \vert} \Phi_j(s,\xi)\, ds 
\right) \quad (t \to \pm\infty). 
\]
Putting 
\[
\Phi_\pm(t,\xi)=
\mathrm{diag}\left(
e^{i\varphi^{\pm}_1(\xi)t},\cdots, 
e^{i\varphi^{\pm}_m(\xi)t} \right),
\]
\[
D_\pm(\xi)=\mathrm{diag}\left(
e^{i\Theta^{\pm}_1(\xi)},\cdots, 
e^{i\Theta^{\pm}_m(\xi)} \right), 
\]
\[
\Psi(t,\xi)=
\mathrm{diag}\left( \Psi_1(t;\xi),\cdots,\Psi_m(t;\xi)\right),
\]
we can write \eqref{phase} as 
\begin{equation}\label{Decomposition}
\Phi(t,\xi)=\Phi_\pm(t,\xi)D_\pm(t,\xi)+\Psi(t,\xi) \qquad \text{with}
\end{equation}
\begin{equation}\label{decay rate}
\Psi(t,\xi)\to 0 \qquad (t\to \pm\infty).
\end{equation}
Plugging this identity into \eqref{EQ:Representation of U} 
from Proposition \ref{prop:Rep}, we have 
\begin{align*}
U(t,x)&=\sum^{m-1}_{j=0} \mathscr{F}^{-1} 
\left[ \mathscr{N}_\pm(\xi)^{-1}
\Phi_\pm(t,\xi)D_\pm(\xi)\pmb{\alpha}^j_\pm(\xi)
\widehat{f}_j (\xi)\right](x)\\ 
&+\sum^{m-1}_{j=0} \mathscr{F}^{-1} 
\left[ \left(\mathscr{N}(t,\xi)^{-1}-\mathscr{N}_\pm(\xi)^{-1}\right)
\Phi_\pm(t,\xi)D_\pm(\xi)\pmb{\alpha}^j_\pm(\xi)
\widehat{f}_j (\xi)\right](x)\\ 
&+\sum^{m-1}_{j=0} \mathscr{F}^{-1} 
\left[ \mathscr{N}_\pm(\xi)^{-1}
\Phi_\pm(t,\xi)D_\pm(\xi)\left(\pmb{a}^j(t,\xi)-\pmb{\alpha}^j_\pm(\xi)\right)
\widehat{f}_j (\xi)\right](x)\\ 
&  +\sum^{m-1}_{j=0} \mathscr{F}^{-1} 
\left[ \mathscr{N}(t,\xi)^{-1}\Psi(t,\xi)\pmb{a}^j(t,\xi) 
\widehat{f}_j (\xi)\right](x), \quad t \gtrless 0.
\end{align*}
It can be readily checked that 
\begin{equation}\label{EQ:V-rep}
V_\pm(t,x)=\sum^{m-1}_{j=0} \mathscr{F}^{-1} 
\left[ \mathscr{N}_\pm(\xi)^{-1}
\Phi_\pm(t,\xi)D_\pm(\xi)\pmb{\alpha}^j_\pm(\xi)
\widehat{f}_j (\xi)\right](x)
\end{equation} 
satisfy the equation \eqref{Limiting Eq}. 
Thus we conclude that 
\[
\Vert U(\cdot,t)-V_\pm(\cdot,t)\Vert_{(L^2(\mathbb{R}^n))^m} \to 0 \quad 
(t \to \pm\infty),
\]
if we use \eqref{decay rate} and the following convergence:
\[
\mathscr{N}(t,\xi)^{-1}\to \mathscr{N}_\pm(\xi)^{-1}, \quad 
\pmb{a}^j(t,\xi)\to \pmb{\alpha}^j_\pm(\xi) \quad 
(t\to\pm\infty).
\]
As a conclusion, $U(x,t)$ is asymptotically free. Moreover, 
the mapping 
$$
\mathscr{W}_\pm^{-1}:U(0)\mapsto 
V_\pm(0)=\sum^{m-1}_{j=0} \mathscr{F}^{-1} 
\left[ \mathscr{N}_\pm(\xi)^{-1}
D_\pm(\xi)\pmb{\alpha}^j_\pm(\xi)
\widehat{f}_j (\xi)\right](x)
$$ 
is bijective and bounded on $(L^2(\Rn))^m$. 
Theorem \ref{thm:Nonfree} (i) is thus proved. \qed 

\section{Proof of Theorem \ref{thm:Nonfree} (ii).} 

We recall from \eqref{EQ:Representation of U} and \eqref{EQ:V-rep} 
that 
\[
\widehat{U}(t,\xi)=\sum_{j=0}^{m-1}\mathscr{N}(t,\xi)^{-1}\Phi(t,\xi)
\pmb{a}^j(t,\xi)\widehat{f}_j(\xi),
\]
\[
\widehat{V}_\pm(t,\xi)=\sum_{j=0}^{m-1}\mathscr{N}_\pm(\xi)^{-1}
\Phi_\pm(t,\xi)D_\pm(\xi)\pmb{\alpha}^j_\pm(\xi)
\widehat{f}_j (\xi).
\]
Thus the proof of 
Theorem \ref{thm:Nonfree} (ii) is reduced to the the next lemma 
provided that $(f_0(x),\ldots,f_{m-1}(x))$ are non-trivial. 
\begin{lem} \label{lem:4.1}
Let $\varphi_j(t,\xi)$ and $\varphi_j^\pm(\xi)$, $j=1,\ldots,m$, 
be the phase functions 
as in \eqref{hyp3} and \eqref{Limiting phase}, respectively. Suppose that 
\begin{equation}
\vert \vartheta_j(t,\xi)-\varphi^{\pm}_j(\xi)t \vert \to +\infty \quad 
(t \to \pm\infty), 
\label{diverge0}
\end{equation}
where $\vartheta_j(t,\xi)=\int^{t}_{0}\varphi_j(s,\xi) \, ds$. 
Let $A_j(t,\xi),B_j(t,\xi)\in C(\R;L^2(\Rn))$ for 
$j=1,\ldots,m-1$, satisfying 
\[
A_j(t,\xi)\to A_j^\pm(\xi), \quad
B_j(t,\xi)\to B_j^\pm(\xi) \quad 
\text{for each $\xi\in \Rn$} \quad 
(t\to\pm\infty),
\]
with some $A_j^\pm(\xi),B_j^\pm(\xi)\in L^2(\Rn)$. 
Then we have 
\begin{multline*}
\left\| \sum^{m}_{j=1} \left\{ 
A_j(t,\xi)e^{i\vartheta_j(t,\xi)}
-B_j(t,\xi)e^{i\varphi^{\pm}_j (\xi)t} \right\} 
\right\|_{L^2(\Rn)} \\
\longrightarrow \left\{
\sum^{m}_{j=1} \left(
\Vert A_j^\pm(\xi)\Vert^{2}_{L^2(\Rn)}
+\Vert B_j^\pm(\xi)\Vert^{2}_{L^2(\Rn)} \right)\right\}^{1/2}
\end{multline*}
as $t\to \pm\infty$. 
\end{lem}
\begin{proof} Putting 
\[
K_\pm(t,\xi)=\sum^{m}_{j=1} \left\{ 
A_j(t,\xi)e^{i\vartheta_j(t,\xi)}
-B_j(t,\xi)e^{i\varphi^{\pm}_j (\xi)t} \right\}, 
\]
we can write 
\begin{equation}
\vert K_\pm(t,\xi)\vert^2=\sum^{m}_{j=1} \left(
\vert A_j(t,\xi)\vert^2+\vert B_j(t,\xi)\vert^2 \right)
+\mathrm{Re} H_\pm(t,\xi), 
\label{spagnolo}
\end{equation}
where 
\begin{multline*}
H_\pm(t,\xi)=2\sum_{j<k} \left\{ 
e^{i\{\vartheta_j(t,\xi)-\vartheta_k(t,\xi)\}}
A_j(t,\xi)\overline{A_k(t,\xi)}
+e^{i\{\varphi^{\pm}_j(\xi)t-\varphi^{\pm}_k(\xi)t\}}
B_j(t,\xi)\overline{B_k(t,\xi)} \right\}\\
-2\sum^{m}_{j=1}
e^{i\{\vartheta_j(t,\xi)-\varphi^{\pm}_j(\xi)t\}}
A_j(t,\xi)\overline{B_j(t,\xi)}
-2\sum_{j<k}
e^{i\{\vartheta_j(t,\xi)-\varphi^{\pm}_k(\xi)t\}}
A_j(t,\xi)\overline{B_k(t,\xi)}.
\end{multline*}
We can check that all the phases in $H_\pm(t,\xi)$ are unbounded in $t$. 
Indeed, it 
follows from \eqref{hyp3} and \eqref{Limiting phase} 
that if $j<k$, then 
\begin{equation}\label{int1} 
\left|\vartheta_j(t,\xi)-\vartheta_k(t,\xi) \right| 
= \Big\vert 
\int^{t}_{0} \left(\varphi_j(s,\xi)-\varphi_k(s,\xi) \right)
\, ds \Big\vert 
\ge d \vert \xi \vert \vert t \vert \to +\infty,
\end{equation}
\begin{equation}
\vert \varphi^{\pm}_j(\xi)t-\varphi^{\pm}_k(\xi)t \vert
\ge d \vert \xi \vert \vert t \vert \to +\infty, 
\label{int2}
\end{equation}
as $t \to \pm\infty$. Since $\varphi_j(t,\xi)\to \varphi^{\pm}_j(\xi)$ 
as $t \to \pm\infty$ by 
Proposition \ref{prop:root}, it follows that for any $\varepsilon>0$ there 
exists a number $T>0$ such that 
\[
\vert \varphi_j(s,\xi)-\varphi^{\pm}_k(\xi)\vert \ge (d-\varepsilon)|\xi|, 
\quad j \ne k, \quad \vert s \vert>T, \quad \xi \ne 0, 
\]
hence, 
\begin{eqnarray}
& & \left\vert \vartheta_j(t,\xi)-\varphi^{\pm}_k(\xi)t\right\vert
= \Big\vert 
\int^{t}_{0} \left(\varphi_j(s,\xi)-\varphi^{\pm}_k(\xi) \right)
\, ds \Big\vert \label{diverge}\\
&\ge& 
\Big\vert
\int^{t}_{T} \left(\varphi_j(s,\xi)-\varphi^{\pm}_k(\xi) \right)
\, ds \Big\vert 
-\Big\vert
\int^{T}_{0} \left(\varphi_j(s,\xi)-\varphi^{\pm}_k(\xi) \right)
\, ds \Big\vert \nonumber \\
&\ge& (d-\varepsilon) |t-T| |\xi| -\Big\vert
\int^{T}_{0} \left(\varphi_j(s,\xi)-\varphi^{\pm}_k(\xi) \right)
\, ds \Big\vert \to +\infty \quad (t \to \pm\infty). \nonumber 
\end{eqnarray}
Now, we note a lemma on oscillatory integrals: \\

\noindent 
{\bf Fact A.} If $\phi(t) \to \pm\infty$ $(t \to \pm\infty)$, then 
\[
\int_{\mathbb{R}^{n}} e^{i\phi(t)\vert \xi \vert}
\Phi(\xi) \, d \xi \to 0, \quad \text{$(t \to \pm\infty)$, \quad 
for $\Phi \in \mathscr{S}(\Rn)$.} 
\]
{} \\

Thus we have, by Fact A and \eqref{diverge0}--\eqref{diverge}, 
\[
\int_{\mathbb{R}^{n}} H_\pm(t,\xi)\, d \xi 
\to 0 \quad (t \to \pm\infty).
\]
In conclusion, we have 
\[
\int_{\mathbb{R}^{n}} 
\vert K_\pm(t,\xi)\vert^2 \, d\xi \to 
\sum^{m-1}_{j=0}\int_{\mathbb{R}^{n}} 
\left(\vert A_j^\pm(\xi) \vert^2+\vert B_j^\pm(\xi)\vert^2 \right) \, d \xi
\quad (t\to \pm\infty).
\]
The proof of Lemma \ref{lem:4.1} is complete. 
\end{proof}

\section{Proof of Theorem \ref{thm:Scattering}.}
\setcounter{equation}{0}
Let $V_\pm=V_\pm(t,x)$ be the solution to the Cauchy problem 
\[
D_tV_\pm=A_\pm(D_x)V_\pm, \quad x\in \Rn, \quad \pm t>0
\]
with Cauchy data 
\[
V_\pm(0,x)={}^T(f_0^\pm(x),\ldots,f_{m-1}^\pm(x)).
\]
Let $\mathscr{N}_\pm(\xi)$ be the diagonaliser of the symbol 
$A_\pm(\xi/|\xi|)$, i.e., 
\[
\mathscr{N}_\pm(\xi)A_\pm(\xi/|\xi|)
=\mathscr{D}_\pm(\xi/|\xi|)\mathscr{N}_\pm(\xi),
\]
where we put 
\[
\mathscr{D}_\pm(\xi)=\mathrm{diag}\left( 
\varphi_1^\pm(\xi),\cdots,\varphi_m^\pm(\xi)
\right).
\]
Denoting 
\[
\Phi_\pm(t,\xi)=\mathrm{diag}\left( 
e^{i\varphi_1^\pm(\xi)t},\cdots,e^{i\varphi_m^\pm(\xi)t}
\right), 
\]
and by $\pmb{e}^0,\ldots,\pmb{e}^{m-1}$ the standard unit vectors 
in $\mathbb{R}^m$, we can write 
\[
V_\pm(t,x)=\sum_{j=0}^{m-1}\mathscr{F}^{-1}
\left[\mathscr{N}_\pm(\xi)^{-1}\Phi_\pm(t,\xi)
\pmb{e}^j \widehat{f}_j^\pm(\xi)\right](x).
\]
In the first step, we will find the asymptotic 
integration for $D_t \widehat{U}(t,\xi)=A(t,\xi)\widehat{U}(t,\xi)$, 
such that $\widehat{U}(t,\xi)$ is asymptotic to 
$\widehat{V}_\pm(t,\xi)$ as $t\to\pm\infty$. 

We prove the following:
\begin{prop} \label{prop:prop5.2}
Assume \eqref{strict hyperbolicity5}--\eqref{strict hyperbolicity6}. 
Let $\mathscr{N}(t,\xi)$ be the diagonaliser of $A(t,\xi/|\xi|)$ from 
Lemma \ref{Diagonalisation}, and put 
\[
\Phi(t,\xi)=\mathrm{diag}\left( 
e^{i\int_0^t\varphi_1(s,\xi)\, ds},\cdots,
e^{i\int_0^t\varphi_m(s,\xi)\, ds}
\right).
\]
Then there exist a 
fundamental matrix $W_\pm(t,\xi)$ of 
$D_t \widehat{U}(t,\xi)=A(t,\xi)\widehat{U}(t,\xi)$ 
such that 
\[
W_\pm(t,\xi)=\mathscr{N}(t,\xi)^{-1}\Phi(t,\xi)(I+R_\pm(t,\xi)) 
\qquad {with}
\]
\[
R_\pm(t,\xi)\to 0 \qquad (t\to\pm\infty).
\]
\end{prop} 
\begin{proof} We prove the case ``$-$", since the case 
``$+$" is the same as ``$-$". The idea of proof comes from 
\cite{Hartman}. 
Let $\varphi_1(t,\xi),\ldots,\varphi_m(t,\xi)$ be the characteristic roots 
of the operator $D_t-A(t,D_x)$. 
We define a matrix $C_-(t,\xi)$ as 
\begin{equation}\label{EQ:C-matrix}
C_-(t,\xi)=\Phi(t,\xi)^{-1}(D_t\mathscr{N}(t,\xi))\mathscr{N}(t,\xi)^{-1}
\Phi(t,\xi).
\end{equation}
Let $\Sigma_-$ be a set of all $\sigma \in \mathbb{R}$ such that 
$\pmb{a}_-^j(t,\xi;\sigma)$, $j=0,\ldots,m-1$, 
are solutions of the problem 
\[
D_t \pmb{a}_-^j(t,\xi;\sigma)=C_-(t,\xi)\pmb{a}_-^j(t,\xi;\sigma), 
\quad 
\pmb{a}_-^j(\sigma,\xi;\sigma)=\pmb{e}^j.
\] 
Hence $\pmb{a}_-^j(t,\xi;\sigma)$ can be written as the Picard series: 
\begin{multline*}
\pmb{a}_-^j(t,\xi;\sigma)=
\left(I+i\int_\sigma^t C_-(\tau_1,\xi) \, d\tau_1 \right. \\
\left. 
+i^2\int_\sigma^t C_-(\tau_1,\xi) \, d\tau_1
\int_\sigma^{\tau_1} C_-(\tau_2,\xi) \, d\tau_2
+\cdots \right)\pmb{e}^j
\end{multline*}
for all $t\in \R$, $\sigma\in \Sigma_-$ and $\xi \in \Rn$. 
Then we can estimate 
\begin{equation}\label{EQ:a-Bound}
\left|\pmb{a}_-^j(t,\xi;\sigma)\right|\le 
e^{\int_\R \|C_-(s,\xi)\|_{L^\infty(\Rn)} \, ds}
\le 
e^{c \int_\R \|\partial_s \mathscr{N}(s,\xi)\|_{L^\infty(\Rn)} \, ds}.
\end{equation}
Put 
$$W_-^j(t,\xi;\sigma)=\mathscr{N}(t,\xi)^{-1}
\Phi(t,\xi)\pmb{a}_-^j(t,\xi;\sigma).
$$ 
Then we see from \eqref{EQ:C-matrix} that 
each $W_-^j(t,\xi;\sigma)$ satisfies the following equation: 
\[
D_tW_-^j(t,\xi;\sigma)=A(t,\xi)W_-^j(t,\xi;\sigma).
\]
It follows from \eqref{EQ:a-Bound} that 
$\pmb{a}_-^j(t,\xi;\sigma)$ exist for all $t \in \mathbb{R}$, 
$j=0,\ldots,m-1$ and 
$\sigma\in \Sigma_-$, and hence, we can estimate 
\[
\| D_t \pmb{a}_-^j(t,\xi;\sigma) \| \le \|C_-(t,\xi)\| 
\| \pmb{a}_-^j(t,\xi;\sigma) \| 
\le c_1e^{c\int_{\R}
\|\partial_\tau \mathscr{N}(\tau,\xi)\|_{L^\infty(\Rn)} \, d \tau}
\|\partial_t \mathscr{N}(t,\xi)\|. 
\]
Since $\partial_t \mathscr{N}(\cdot,\xi)\in L^1(\R)$, these estimates
imply that, for $|t|\le|\sigma|$, 
$\sigma\in\Sigma_-$,
\begin{align} 
\|\pmb{a}_-^j(t,\xi;\sigma)-\pmb{e}^j\| 
=&\left\| \int^{t}_{\sigma} \partial_t \pmb{a}_-(\tau,\xi;\sigma) \, d\tau 
\right\| \label{EQ:Conv} \\ 
\le& c_1 e^{c\int_{\R}
\|\partial_\tau \mathscr{N}(\tau,\xi)\|_{L^\infty(\Rn)} \, d \tau}
\int_{|t|}^{|\sigma|} \|\partial_s \mathscr{N}(s,\xi)\| \, ds. \nonumber
\end{align}
In particular, the family 
$\{ \pmb{a}_-^j(t,\xi;\sigma) \}_{\sigma \in \Sigma_-}$ 
is uniformly bounded in $\sigma,\xi$ 
and equicontinuous on every bounded $t$-interval. Hence 
there exists a sequence $\{ \sigma_\ell\}^{\infty}_{\ell=1}$ 
of $\Sigma_-$ such that 
\[
|\sigma_1|<|\sigma_2|<\cdots, \quad |\sigma_\ell| \to \infty 
\quad (\ell \to \infty), 
\]
and the limits 
\[
\pmb{a}_-^j(t,\xi)=\lim_{\ell\to \infty} \pmb{a}_-^j(t,\xi;\sigma_\ell) 
\]
exist uniformly in $\xi$ 
on every bounded $t$-interval. Moreover, the limits 
\[
W_-^j(t,\xi)
=\lim_{\ell\to \infty}W_-^j(t,\xi;\sigma_\ell)
=\mathscr{N}(t,\xi)^{-1}\Phi(t,\xi)\pmb{a}_-^j(t,\xi)
\]
also exist. Hence $\pmb{a}_-^j(t,\xi)$ are the 
solutions of $D_t\pmb{a}_-^j(t,\xi)=C_-(t,\xi)\pmb{a}_-^j(t,\xi)$, 
and $W_-^j(t,\xi)$ are the solutions 
of $D_tW_-^j(t,\xi)=A(t,\xi)W_-^j(t,\xi)$. 
Putting $\sigma=\sigma_\ell$ in \eqref{EQ:Conv}, and 
letting $\ell \to \infty$, with $t$ 
fixed, we see that 
\begin{equation}
\| \pmb{a}_-^j(t,\xi)-\pmb{e}^j \| 
\le c_1 e^{c\int_{\R}
\|\partial_\tau \mathscr{N}(\tau,\xi)\|_{L^\infty(\Rn)} \, d \tau}
\int^{t}_{-\infty} \|\partial_s \mathscr{N}(s,\xi)\| \, ds
\label{order2} 
\end{equation}
for all $t\in \mathbb{R}$. The uniqueness of each $\pmb{a}_-^j(t,\xi)$ 
is obvious. 

Now we can write, by putting $\pmb{r}_-^j(t,\xi)=\pmb{a}^j(t,\xi)
-\pmb{e}^j$, 
\begin{equation}
W_-^j(t,\xi)=\mathscr{N}(t,\xi)^{-1}\Phi(t,\xi)
(\pmb{e}^j+\pmb{r}_-^j(t,\xi)),
\label{combination}
\end{equation}
where $\pmb{r}_-^j(t,\xi)$ is uniform in $\xi$ and satisfies 
\begin{equation}\label{EQ:r-decay}
\pmb{r}_-^j(t,\xi)\to \pmb{0} \qquad (t\to-\infty)
\end{equation}
on account of \eqref{order2}. 
It remains to prove that 
$W_-(t,\xi)\equiv (W_-^1(t,\xi),\ldots,W_-^{m-1}(t,\xi))$ is the 
fundamental matrix for $D_t \widehat{V}=A(t,\xi)\widehat{V}$. 
Taking the determinant of $W_-(t,\xi)$ in \eqref{combination}, we have, 
by using \eqref{EQ:r-decay}, 
$\mathscr{N}(t,\xi)^{-1} \to \mathscr{N}_-(\xi)^{-1}$ 
$(t\to-\infty)$ and $|\mathrm{det}\Phi(t,\xi)|=1$, that 
\[
|\mathrm{det}W_-(t,\xi)|
\longrightarrow |\mathrm{det}\mathscr{N}_-(\xi)|^{-1}
\ne0 \quad (t\to-\infty).
\]
Hence there exists a number 
$t_0 \le 0$ such that 
\[
\mathrm{det}W_-(t,\xi)\ne0 \quad \text{for all $t \le t_0$.}
\]
Since $W_-(t,\xi)$ satisfies $D_t W_-(t,\xi)=A(t,\xi)W_-(t,\xi)$, 
it follows from the Abel-Jacobi formula that 
\[
\mathrm{det}W_-(t,\xi)=\mathrm{det}W_-(t_0,\xi)
\exp \int^{t}_{t_0} \mathrm{tr}A(s,\xi) \, ds 
\ne 0 
\]
for all $\xi \in \mathbb{R}^n \setminus 0$ and all $t \in \mathbb{R}$. 
It means that $W_-(t,\xi)$ is the fundamental matrix for 
$D_t \widehat{U}=A(t,\xi)\widehat{U}$. This completes 
the proof of Proposition \ref{prop:prop5.2}. 
\end{proof} 

\begin{proof}[Completion of the proof of Theorem \ref{thm:Scattering}]
Now we are in a position to prove our theorem. Namely, 
we will find a solution $U(t,x)$ of 
$D_t U=A(t,D_x)U$ such that 
\begin{equation}\label{aim}
\|V_\pm(t,\cdot)-U(t,\cdot) \|_{(L^2(\mathbb{R}^n))^m} \to 0 
\quad (t \to \pm\infty). 
\end{equation}
Now going back to \eqref{phase}--\eqref{decay rate} in the proof of 
Theorem \ref{thm:Nonfree} (i), we have 
\begin{align*}
W_\pm(t,\xi)=& 
\mathscr{N}_\pm(\xi)^{-1}\Phi_\pm(t,\xi)D_\pm(\xi)
 +(\mathscr{N}(t,\xi)^{-1}-\mathscr{N}_\pm(\xi)^{-1})
\Phi_\pm(t,\xi)D_\pm(\xi)\\
&+\mathscr{N}(t,\xi)^{-1}\Psi(t,\xi).
\end{align*}
Thus putting 
$$\widehat{U}(t,\xi)=\sum_{j=0}^{m-1}W_\pm(t,\xi)D_\pm(\xi)^{-1}
\pmb{e}^j \widehat{f}_j^\pm(\xi), 
$$ 
we can decompose $\widehat{U}(t,\xi)$ into three terms:
\begin{align*}
\widehat{U}(t,\xi)=&\widehat{V}_\pm(t,\xi)
+\sum_{j=0}^{m-1}(\mathscr{N}(t,\xi)^{-1}-\mathscr{N}_\pm(\xi)^{-1})
\Phi_\pm(t,\xi)\pmb{e}^j \widehat{f}_j^\pm(\xi)\\
&+\sum_{j=0}^{m-1}\mathscr{N}(t,\xi)^{-1}\Psi(t,\xi)D_\pm(\xi)^{-1}
\pmb{e}^j \widehat{f}_j^\pm(\xi).
\end{align*}
It can be readily checked that this $\widehat{U}(t,\xi)$ satisfies 
$D_t\widehat{U}(t,\xi)=A(t,\xi)\widehat{U}(t,\xi)$. 
Since $\mathscr{N}(t,\xi)^{-1}-\mathscr{N}_\pm(\xi)^{-1}\to 0$ and 
$\Psi(t,\xi)\to 0$ as $t\to \pm \infty$, we conclude from the 
Plancherel theorem 
that \eqref{aim} is true. 
Moreover, the mapping 
\[
\mathscr{W}_\pm:V_\pm(0)\mapsto 
U(0)=\sum_{j=0}^{m-1} \left[W_\pm(0,\xi)D_\pm(\xi)^{-1}
\pmb{e}^j \widehat{f}_j^\pm(\xi)\right]
\]
is bijective and bounded on $(L^2(\Rn))^m$. 
The proof of Theorem \ref{thm:Scattering} is finished. 
\end{proof} 

\noindent 
{\bf Acknowledgements.} The proof of Lemma \ref{lem:4.1} 
is based on the private communication with 
Professors Sergio Spagnolo and Marina Ghisi when the first 
author visited the University of Pisa. 
The authors would like to express their 
sincere gratitude to them.



\begin{thebibliography}{99} 

 \bibitem{Arosio}
A.~Arosio, {\em Asymptotic behaviour as $t\rightarrow +\infty $ of the 
solutions of linear hyperbolic equations with coefficients discontinuous in 
time $($on a bounded domain$)$}, J. Differential Equations {\bf 39} (1981), 
291--309. 

%
\bibitem{Ascoli}G. Ascoli, 
{\em Sulla forma asintotica degli integrali dell'equazione differenziale 
$y^{\prime \prime}+A(x)y=0$ in un caso notevole di stabilit\`a}, Univ. Nac. 
Tucum\'an, Revista A. {\bf 2} (1941), 131--140. 


%
\bibitem{Coddington}
   E.A. Coddington and N. Levinson, 
``Theory of differential equations'', 
New York, McGraw-Hill, 1955. 

%
 \bibitem{Hartman} 
  P.~Hartman, ``Ordinary differential equations'', SIAM, 2nd edition, 2002. 


%
\bibitem{Matsuyama0}T.~Matsuyama, {\em Asymptotic behaviour for wave equations 
with time-dependent coefficients}, Annali Universit\`a di Ferrara 
Sec. VII  - Sci. Math., {\bf 52} (2) (2006), 383--393. 

%
\bibitem{Matsuyama22}T.~Matsuyama, 
{\em Asymptotic profiles for Kirchhoff equation}, Rend. Lincei 
Mat. Appl. {\bf 17} (2006), 377--395.


\bibitem{MR2}T.~Matsuyama and M. Ruzhansky, {\em Asymptotic integration and dispersion 
for hyperbolic equations},
Adv. Differential Equations {\bf 15} (2010), 721--756.  

\bibitem{MR3}T.~Matsuyama and M. Ruzhansky, {\em Dispersion for 
hyperbolic systems, with applications to Kirchhoff systems}, preprint. 

%
\bibitem{Mizohata}S.~Mizohata, 
``The theory of partial differential equations'', 
Cambridge Univ. Press, 
1973. 

\bibitem{RS}
M.~Ruzhansky and J.~Smith, ``Dispersive and Strichartz estimates for 
hyperbolic equations with constant coefficients'', MSJ Memoirs, 22, 
Mathematical Society of Japan, Tokyo, 2010. 

%
 \bibitem{Wintner}A. Wintner, 
 {\em Asymptotic integrations of adiabatic oscillator}, Amer.\ J.\ Math.\ 
 {\bf 69} (1947), 251--272.


\end{thebibliography}
\end{document}